\documentclass[12pt]{amsart}

\usepackage{amssymb, amsmath, latexsym, amscd}

\newtheorem*{thmA}{Theorem A}
\newtheorem*{thmB}{Theorem B}
\newtheorem*{propC}{Proposition C}
\newtheorem*{propD}{Proposition D}
\newtheorem*{thmE}{Theorem E}
\newtheorem*{propF}{Proposition F}
\newtheorem*{propG}{Proposition G}

\newtheorem{definition}{Definition}
\newtheorem{corollary}{Corollary}
\newtheorem{lemma}{Lemma}
\newtheorem{proposition}{Proposition} 

\newcommand{\GL}{\text{\rm{GL}}}
\newcommand{\PGL}{\text{\rm{PGL}}}

\newcommand{\Z}{{\mathbb Z}}
\newcommand{\Aut}{\text{\rm Aut}}
\newcommand{\Out}{\text{\rm Out}}

\def\iso{\cong}

\def\l{\lambda}
\def\autn{\Aut(F_n)}
\def\e{\varepsilon}

\begin{document}

\title{Homomorphisms from automorphism groups of free groups} 

\thanks{The research of
the first author was supported
by an EPSRC Advanced Fellowship.
The research of the second author is supported in part by NSF grant
DMS-9307313.}

\author[Bridson]{Martin R.~Bridson}
\address{MRB:
Mathematics Department,
180 Queen's Gate,
London SW7 2BZ,
U.K.}
\email{m.bridson@ic.ac.uk}

\author[Vogtmann]{Karen Vogtmann}
\address{KV: Mathematics Department, 555 Malott Hall, Cornell University,
Ithaca, NY 14850, USA}
\email{vogtmann@math.cornell.edu}

\subjclass{Primary: 20F65, 20F28}

\date{\today}

\keywords{automorphism groups of free groups, super-rigidity}

\begin{abstract}{The automorphism group of a finitely generated free group
is the normal closure
of a single element of order $2$.

If $m<n$ then a  homomorphism $\Aut(F_n)\to \Aut(F_m)$ can have
image of  cardinality at most $2$. More generally, this is true
of homomorphisms from $\Aut(F_n)$ to any group  that
does not contain an isomorphic image of the symmetric group
$S_{n+1}$.
Strong restrictions are also obtained on maps
to groups which do
not contain a copy of $W_n=(\Z/2)^n\rtimes S_{n}$, or of $\Z^{n-1}$.

These results place constraints on how $\Aut(F_n)$ can
act. For example, if
   $n\ge 3$ any action of $\Aut(F_n)$  on the
circle (by
homeomorphisms) factors through $\text{\rm{det}} : \Aut(F_n)\to\Z_2$.}
\end{abstract}

\maketitle

\def\QED{\end{proof}}


\section{Introduction}

In recent articles we began to explore the extent to
which the well-known analogies between lattices
in semisimple Lie groups and automorphism groups
of free groups can be extended to cover various aspects of rigidity.

For example,  in \cite{BV} we proved that all automorphisms
of $\Aut(F_n)$ and $\Out(F_n)$ are inner if $n\ge 3$. In
the direction of super-rigidity,
it follows from the main theorem of \cite{BFH}
that if $\Gamma$ is an irreducible, non-uniform
lattice in a higher rank semisimple group, then
any homomorphism from $\Gamma$ to $\Aut(F_n)$
or $\Out(F_n)$ has finite
image --- see \cite{BF}. In the classical setting, Margulis super-rigidity
tells one that if there is no homomorphism from one
semi-simple group $G_1$ to another $G_2$, then any map
from a lattice in $G_1$ to $G_2$ must have finite image. An
   example of this phenomenon is the fact
that if $m<n$ then any map from $\GL (n,\Z)$ to
$\GL (m,\Z)$ has image $\Z_2$ or $\{1\}$. In the
present article we establish the
analogous result for automorphism groups of free groups.

\begin{thmA}  If $n\geq 3$ and $n > m$, then every homomorphism
$\Aut(F_n)\to \Aut(F_m)$ is either trivial or has image of
order $2$. Likewise for maps $\Aut(F_n)\to \Out(F_m)$ and
$\Aut(F_n)\to \GL (m,\Z)$.
\end{thmA}

In section 3 we shall see that Theorem A is a special case of:

\begin{thmB} Let $n\geq 3$. If a group $G$ does not contain a
copy of the symmetric group $S_{n+1}$, then the image of any homomorphism
$\Aut(F_n)\to G$ has cardinality at most $2$.
\end{thmB}

Note that there is only one surjection $\Aut(F_n)\to\Z_2$,
namely the composition of the determinant map $\GL(n,\Z)\to\Z_2$
and the map $\Aut(F_n)\to\GL(n,\Z)$
describing the action of $\Aut(F_n)$ on the abelianization
of $F_n$. We shall denote this surjection  $\text{\rm{det}}:
\Aut(F_n)\to\Z_2$.

Theorem B is proved by examining the pattern
of finite  subgroups in  $\Aut(F_n)$. Another result in
this direction concerns the largest finite subgroup $W_n\subset
\Aut(F_n)$ (which is unique up to conjugacy).

\begin{propC}  Let $n\geq 3$, and let $\phi\colon\Aut(F_n)\to G$ be
a homomorphism to
a group
$G$ that does not contain a  copy of $W_n=(\Z_2)^n\rtimes S_n$.  Then
either $\hbox{\rm{Im}}(\phi)
\iso \PGL (n,\Z)$ or else $\hbox{\rm{Im}}(\phi)$ is finite.
\end{propC}

Our proof of the above results proceeds via the following observation:

\begin{propD} 
 If $n\geq 3$, then $\Aut(F_n)$ (and hence each of its
quotients) is the normal closure of a single element of
order $2$.
\end{propD}

Whenever one is able to control the nature of the quotients
that a group admits, one immediately obtains constraints on
the type of actions that it admits. For example, because
the  finite subgroups of the homeomorphism group
of the circle are all cyclic or dihedral, Theorem B implies:

\begin{thmE}
 If $n\geq 3$, then any action of $\Aut(F_n)$  on the
circle (by
homeomorphisms) factors through $\text{\rm{det}} : \Aut(F_n)\to\Z_2$.
\end{thmE}

B. Farb and J. Franks \cite{FF} give a different proof of this fact for
$C^2$-actions in the case $n\geq 6$.

Our results limiting the actions of $\Aut(F_n)$
do not compare favourably with results
from the classical setting   in that they
rely heavily on the nature of the torsion in $\Aut(F_n)$
and therefore do  not extend to subgroups of finite index (which
may be torsion-free).  One might be tempted to
conjecture   that if $\Gamma\subset\Aut(F_n)$ has finite
index  and if $m<n$ then every homomorphism from $\Gamma$
to $\Aut(F_m)$
has finite image. However this fails for $n=3$,
because   $\Aut(F_3)$ has a torsion-free
subgroup of finite index which maps onto $\Z$. The situtation
for $n > 3$ is far from clear.

In the case of $\GL (n,\Z)$ it is easy to extend results to
finite index subgroups, because one has a very concrete description
of these subgroups in terms of congruence subgroups.
Part of the difficulty for $\Aut(F_n)$ is that the following naive
analogues $Q_{n,m}$ of the congruence quotients
$\GL (n,\Z/m\Z)$ can be  infinite.

Fix a basis $\{a_1,\dots,a_n\}$ for $F_n$ and let
$\lambda_{ij}$ be the automorphism of $F_n$ that sends $a_i$ to
$a_ja_i$ and
fixes $a_k$ for $k\neq i$.  Let $Q_{n,m}$ be the quotient of
$\Aut(F_n)$ defined by
setting $\lambda_{ij}^m=1$ for all $i$ and $j$.

\begin{propF}
$Q_{n,m}$ surjects onto a group that contains a copy of
the free Burnside group of exponent $m$
on $n-1$ generators.
\end{propF}

The groups $Q_{n,m}$ enjoy the following universal property
with respect to maps from $\Aut(F_n)$ to groups without
large-rank abelian subgroups (see section 5).

\begin{propG}
If  $G$ is a group that  does not contain $\Z^{n-1}$,
then any homomorphism $\Aut(F_n) \to G$
factors through
$Q_{n,m}$ for some $m$.
\end{propG}

\section{Some elements and relations in $\Aut(F_n)$}

We assume that $n\ge 3$  and we fix a basis
$\{a_1,\dots,a_n\}$ for the free group
$F_n$.  Associated to this choice of basis one has the
finite subgroup isomorphic to $\Z_2^n$ generated by the involutions
$\e_i$, where $\e_i\colon
   a_i\mapsto a_i^{-1}$ and $a_j\mapsto a_j$ for $j\neq i$.
   The permutations of
our fixed basis form a copy of the
symmetric group $S_n\subset \Aut(F_n)$. We shall write elements
of this symmetric group as products of cycles, for example
$(1\ 2)$ will denote the automorphism that interchanges $a_1$
and $a_2$ and leaves the remaining $a_i$ fixed.

Note that $\e_{\sigma(i)}=\sigma\e_i\sigma^{-1}$
for each permutation $\sigma\in S_n$.

The $\e_i$ and $S_n$
together generate a subgroup
$W_n\cong (\Z_2)^n\rtimes S_n$, which is the unique largest finite subgroup
of $\Aut(F_n)$ for $n\geq 4$.

Associated to our basis we also have the left Nielsen transformations:
$\lambda_{ij}$ sends $a_i$ to
$a_ja_i$ and
fixes $a_k$ for $k\neq i$. One also has the
corresponding right Nielsen moves $\rho_{ij}: a_i\mapsto a_ia_j$.
All Nielsen moves lie in
the index 2 subgroup $\Aut^+(F_n)$  which is
the inverse image of $\text{\rm{SL}}(n,\Z)$ under the natural map
$\Aut(F_n)\to \GL (n,Z)$. It is well-known that they
generate $\Aut^+(F_n)$ (see, e.g., \cite{Gersten}).

If $\sigma$ is a permutation that sends $(i,j)$ to
$(k,l)$, then $\sigma$ conjugates $\l_{ij}$ to
$\l_{kl}$. And   conjugation by $\e_i\e_j$ sends
$\lambda_{ij}$ to $\rho_{ij}$.
Thus all Nielsen moves with respect to a fixed
basis are conjugate. And since $\Aut(F_n)$ acts transitively
on bases of $F_n$, the Nielsen moves associated to different
bases are also all conjugate in $\Aut(F_n)$. In particular
we have:

\begin{lemma}\label{A} If one  Nielsen
move lies in the kernel of a homomorphism
$\phi:\Aut(F_n)\to G$, then $\phi$ factors through
 $\text{\rm{det}}:\Aut(F_n)\to\Z_2$. 
\end{lemma}

In order to write relations in $\Aut(F_n)$ in a
standard and convenient manner, it is best to think
of $\Aut(F_n)$ acting {\it on the right} on $F_n$,
and we shall adopt this  convention (which is standard
in texts on combinatorial group theory). Our
commutator convention is $[a,b]=aba^{-1}b^{-1}$.
Thus
$[\alpha,\beta]$ means {\it apply the automorphism $\alpha$,
then $\beta$, then $\alpha^{-1}$, then $\beta^{-1}$}.

We shall make frequent use of the relations
$$\ [\lambda_{i j}, \lambda_{j k}] = \lambda_{i k}\
\ \hbox{and}\ \  [\lambda_{i j}, \lambda_{k j}] = 1.
$$
for  $i,j$ and $k$  distinct.

\begin{proposition}\label{Sn}
 Let $n\geq 3$ and let $\phi:
\Aut(F_n)\to G$ be a homomorphism.  If $\phi$ is not injective on $S_n$,
then
$\phi$ is either trivial or has image $\Z_2$.
\end{proposition}

\begin{proof} Let $K$ be the kernel of $\phi|_{S_n}$.  Since $n\geq
3$ and $K\neq S_n$, we must
have $K=A_n$ or, if $n=4$, possibly $K=\Z/2\times \Z/2$.

If $K=A_n$, then all 3-cycles $(i\,j\,k)$ are in the kernel, so the
relations
$$(i\, j\, k) \lambda_{jk} (i\, j\, k)^{-1}= \lambda_{ij} \hbox{
and } \ [\lambda_{ij}, \lambda_{jk}] = \lambda_{ik}$$ show
that $\lambda_{ik}$ maps trivially under $\phi$, and hence the whole of
$\Aut^+(F_n)$
maps trivially (Lemma \ref{A}).

If  $n=4$ and $K=\Z/2\times\Z/2$, then  all products of two disjoint
2-cycles in
$S_4$ are mapped trivially, so  the identity
$$(j\, k)(i\, l)\  \lambda_{j k}\  (j\, k)(i\, l) = \lambda_{k j}$$
gives    $\phi(\lambda_{j k})=\phi(\lambda_{k j})$.  Applying $\phi$
to the relation
$$[\lambda_{i j}, \lambda_{j k}] = \lambda_{i k}$$ gives
$$\phi(\lambda_{i k})=\phi([\lambda_{i j},
\lambda_{j k}])=\phi([\lambda_{i j}, \lambda_{k j}])=1,$$
and as before all of $\Aut^+(F_n)$ is mapped trivially.
 \QED

Let $\iota = (1\, 2)$.

\begin{proposition} Let $G$ be a group and let $\phi:
\Aut(F_n)\to G$ be a homomorphism. The image of
$\phi$ is trivial if and only if $\phi(\iota)=1$.
\end{proposition}

\begin{proof} If $\phi(\iota) =1$, then since
$S_n$ is generated by conjugates of $\iota$, the
whole of $S_n$ has trivial image. The previous
proposition then shows that $\phi(\Aut^+(F_n))$
is trivial, and since $\iota\notin \Aut^+(F_n)$,
in fact $\phi(\Aut(F_n))$ must be trivial.
\QED

\begin{corollary} $\autn$ is the normal closure of the transposition
$(1\,2)$.
\end{corollary}

\section{Maps to groups without large symmetric groups}

We shall exploit the following well-known facts.

\begin{lemma} $\GL (n-1,\Z)$ does not contain a copy
of $S_{n+1}$.
\end{lemma}
\begin{proof} See, e.g., \cite{FH}, Exercise 4.14.\QED

\begin{lemma}\label{sigm} 
$\autn$ contains a symmetric group
$\Sigma\cong S_{n+1}$ that intersects $W_n=(\Z_2)^n\rtimes S_n$
in the visible copy of $S_n$.
\end{lemma}

\begin{proof}
 Finite subgroups of $\Aut(F_n)$ correspond to
vertex stabilizers in Outer Space \cite{CV}.
The subgroup $\Sigma$ fixes the two vertices of the marked
graph that has $(n+1)$ directed edges with common source and sink,
 $n$ of which are
labelled $a_1,\dots, a_n$.
\QED

\begin{thmB}\label{Snplusone}
     Let $G$ be a group, and suppose $n\geq 3$. If $G$ does not contain a
copy
of $S_{n+1}$ then every homomorphism $\phi:\autn \to G$
has image of cardinality at most $2$.
\end{thmB}

\begin{proof}
Let $\Sigma$ be as in the previous lemma, and let $K$ be the kernel
of $\phi|_\Sigma$.
If $K=A_{n+1}$, then $K\cap S_n=A_n$, so $\phi$ is not injective on
$S_n$ and we are done by
Proposition \ref{Sn}.  There is an additional possibility when $n=3$, namely
$K\cong \Z/2\times\Z/2$, generated by
the order 2 automorphisms $\alpha$ and $\beta$:

$$
{\alpha\colon\begin{cases}  a_1\mapsto a_1^{-1}&\cr
            a_2\mapsto a_3a_1^{-1}\cr
            a_3\mapsto a_2a_1^{-1}
            \end{cases}}
\qquad\qquad{\beta\colon\begin{cases}  a_1\mapsto a_3a_2^{-1}\cr
            a_2\mapsto a_2^{-1}\cr
            a_3\mapsto a_1a_2^{-1}\end{cases}}$$

The relation
$$\beta\lambda_{12}\beta=\lambda_{32}^{-1}$$
gives $\phi(\lambda_{12})= \phi(\lambda_{32}^{-1})$. But expanding
the relation
$$\lambda_{12} = [\lambda_{13},\lambda_{32}]$$
gives
$$
\lambda_{12} = (\lambda_{13} \lambda_{32} \lambda_{13}^{-1})\,
\lambda_{32}^{-1},
$$
so applying $\phi$ we deduce
$\phi(\lambda_{13} \lambda_{32} \lambda_{13}^{-1})=1$,
and $\phi(\lambda_{32})=1$. Thus, in the light of Lemma \ref{A},
the image of $\phi$ is at most $\Z_2$.
\QED

\begin{thmA} If $n\geq 3$ and $m<n$, then every homomorphism $\autn
\to
\Aut(F_m)$ has image of cardinality at most $2$.
\end{thmA}

\begin{proof}
The kernel of the natural map $\Aut(F_m)\to \GL (m,\Z)$
is torsion-free \cite{BT}, so the first of the above lemmas tells
us that $\Aut(F_m)$ does not contain a copy of $S_{n+1}$.
\QED

Exactly the same argument applies to maps from $\Aut (F_n)$ to
$\Out(F_m)$ or $\GL (m,\Z)$.

\section{Maps to groups without copies of $W_n$}

\begin{lemma}\label{KinN}  
Consider a semi-direct product
$N\rtimes S_n$ where $N$ is abelian and one of the following conditions
holds: $n\ge 5$, or $n=4$ and $N$ does not  contain an
$S_n$-invariant subgroup isomorphic to $\Z_2\times\Z_2$,
or $n=3$ and $N$ does not  contain an $S_n$-invariant
subgroup of order 3.
If
$\phi\colon N\rtimes S_n\to Q$ is a homomorphism
with non-trivial kernel $K$, then either $K\cap S_n$
or $N\cap K$ is non-trivial.
\end{lemma}

\begin{proof} If both  $K\cap S_n$
and $N\cap K$  are trivial, then since $\phi$ is assumed to
be non-injective, $\phi(N)\cong N$ and $\phi(S_n)\cong S_n$
must intersect non-trivially. Since $N$ is normal in  $N\rtimes S_n$,
this intersection is normal in $\phi(S_n)$. If $n\ge 5$, then
$S_n$ has no non-trivial normal abelian subgroups, so we have
a contradiction. In $S_4$  one has the possibility that
the intersection could be the Klein 4-group, and in $S_3$ it
could be $A_3\cong \Z_3$.

But the intersection is of
the form $\phi(N_0)\cong N_0$, where $N_0\subset N$ is an $S_n$-invariant
subgroup of $N$. And we have imposed   hypotheses to exclude
the possible existence of non-trivial $N_0$ in the cases $n=3$
and $n=4$. (We leave the reader to consider the easy case $n=2$.)
\QED

\begin{propC} Let $n\geq 3$. If
$\phi\colon\Aut(F_n)\to G$ is a homomorphism to
a group
$G$ that does not contain a  copy of $(\Z_2)^n\rtimes S_n$, then
the image of $\phi$ is either isomorphic to $\PGL (n,\Z)$ or
else it  is finite.
\end{propC}

\begin{proof}
  Since $G$ does not contain a subgroup isomorphic to
$W_n=(\Z_2)^n\rtimes S_n$, the map
$\phi|_{W_n}$ has non-trivial kernel $K$.   If
     $\phi|_{S_n}$ is not injective, then the image of $\phi$ has order
at most $2$, by Proposition \ref{Sn}.

Now assume $\phi$ is injective on $S_n$.
The only non-trivial  subgroups of
$N:=(\Z_2)^n\subset W_n$  that are $S_n$-invariant are  $N$,
its centre
$\langle z\rangle$ ($z=\e_1\e_2\ldots\e_n$) and
$H=\langle
\e_i\e_j \mid i\neq j\rangle$. These are all 2-groups, of course,
but in the case $n=4$ none of them is isomorphic to $\Z_2\times\Z_2$.
Thus we are in the situation of the lemma, and deduce that
$K:=\ker\phi|_N$ is non-trivial. Since $K$ is normal in $W_n$,
it must be one of the subgroups of $N$ listed above.

If $K= N$ or $\langle z\rangle$, or if $n$ is even and $K=H$, then
$\phi(z)=1$, and the identity
$z\lambda_{ij}z=\rho_{ij}$ implies that
$\phi(\lambda_{ij}) =\phi(\rho_{ij})$ for all $i,j$. Since adding the
relations $\lambda_{ij}=\rho_{ij}$ to the standard
presentation for $Aut(F_n)$ gives a presentation for $\GL (n,\Z)$,
this implies that $\phi$ factors through $\GL (n,\Z)$.  Since $\phi(z)$
must be contained in the centre
of $\GL (n,\Z)$, in fact $\phi$ factors  
through $P\GL (n,\Z)$.  And since all normal subgroups of $P\GL (n,\Z)$ have
finite index \cite{B}, the
image of $\phi$ is either finite or isomorphic to
$P\GL (n,\Z)$.

If $n$ is odd and $K=H$, then
$K$ does not contain $z$. But  $\phi(\e_i\e_j)=1$ for all $i,j$, so
the relations
$\e_i\e_j\lambda_{ij}\e_j\e_i=\rho_{ij}$
again imply that $\phi(\lambda_{ij})=\phi(\rho_{ij})$ for all $i,j$.
Therefore
$\phi$ factors through $\GL (n,\Z)$.
Let $\bar \phi\colon \GL (n,\Z)\to G$ be the induced map. Since
$\phi(\e_1\e_2)=1$,
the commutator $[\lambda_{12},\e_1\e_2]$
is in the kernel of $\phi$.  The image  of $[\lambda_{12},\e_1\e_2]$ in
$\GL (n,\Z)$ is
the matrix $I + 2E_{21}$,
which has infinite order in $\GL (n,\Z)$.
Since
    $\GL (n,\Z)$
has no infinite
normal subgroups of infinite index, the image of $\phi$ must be finite.
 \QED

\section{Maps to groups without copies of $\Z^{n-1}$}

We now consider maps of $\autn$ to groups $G$ which do not contain
large free abelian subgroups. We shall see that the following family
of quotients of $\Aut(F_n)$ play a distinguished role in the
study of such maps.

\begin{definition}   Let $Q_{n,m}$ be the quotient of $\Aut(F_n)$ by the
relations $\lambda_{ij}^m=1$ for all $i\neq j$.
\end{definition}

\begin{proposition} Let $n\geq 3$. If $G$ is a group that does not
contain a free
abelian group of rank $n-1$, then every homomorphism
$\phi:\Aut(F_n)\to G$ factors through $Q_{n,m}$ for some $m\ge 1$.
\end{proposition}

\begin{proof} Let $S_{n-1}$ be the symmetric group which permutes the
basis elements
$\{a_2,\ldots,a_n\}$ of $F_n$,  let $N\cong \Z^{n-1}$ be the subgroup of
$\autn$  generated by the elements $\lambda_{i 1}$, $i=2,\ldots,n$,
and let $L=N\rtimes
S_{n-1}$.

Both $N$ and the kernel of $\phi$ are normalized by
$S_{n-1}$, and hence so is their intersection, which we denote $K$.
(This is non-trivial by hypothesis.)
The only
$S_{n-1}$-invariant subgroups of $N\cong \Z^{n-1} $ are the lattices
$mN$, for $m\in \Z$, the diagonal subgroup
$D=\langle(\lambda_{21},\ldots,\lambda_{n1})\rangle\cong \Z$, the
complementary hyperplane
$H=\{(\lambda_{21}^{m_1},\ldots,\lambda_{n1}^{m_n})|\sum
m_i=0\}\cong \Z^{n-2}$,  and  the intersections $mN\cap D$ and $mN\cap H$;
thus these are the only possibilities for  $K$.

\noindent{\bf Case 1.}   $K=mN$.

This case is trivial;  $\phi$ sends all
$\lambda_{ij}^m$ to 1, so factors through $Q_{m,n}$.  

\noindent{\bf Case 2.}   $K= mN\cap H$. 

Since
$mN\cap H$ contains the elements
$\lambda_{i1}^m\lambda_{j1}^{-m}$ for $i\neq j$, we have
$\phi(\lambda_{i1}^m)=\phi(\lambda_{j1}^m)$ for all $i,j$.

A simple calculation shows that 
$[\l_{12},\l_{23}^m]=\l_{13}^m$. By  taking the
image of this relation  under $\phi$ and making the substitution
$\phi(\l_{13}^m)=\phi(\l_{23}^m)$ we get
 $\phi(\l_{12})\phi(\l_{13})^m\phi(\l_{12})^{-1}=
\phi(\l_{13})^{2m}$.
Similarly, from  $[\l_{13},\l_{32}^m]=\l_{12}^m$
and $\phi(\l_{12}^m)=\phi(\l_{32}^m)$, we get 
 $\phi(\l_{13})\phi(\l_{12})^m\phi(\l_{13})^{-1} =
\phi(\l_{12})^{2m}$.

Set $a=\phi(\lambda_{12})$ and $b = \phi(\lambda_{13})^m$. Then
$aba^{-1}=b^2$ and
$ba^mb^{-1} = a^{2^mm}$.  These relations imply that  
$a=\phi(\lambda_{12})$ has finite order. Indeed, in any group,
if $aba^{-1}=b^2$ and $ba^rb^{-1}=a^s$ with 
$r\neq s$,
then  $b$ (and hence $a$) has finite order,
because if we conjugate $b$ by the left-hand side of 
the second relation then the first relation tells us that
the result is $b^{2^r}$, whereas if we  conjugate by
the right-hand side of the second relation we
obtain $b^{2^s}$. 

Since all $\lambda_{ij}$ are conjugate in $\Aut(F_n)$, we conclude that
$\phi(\lambda_{ij}^{m'})=1$ 
for all $i$ and $j$ and some $m'$, i.e. $\phi$ factors through $Q_{n,m'}$.

\noindent{\bf Case 3.}   $K=mD=mN\cap D$.

The diagonal $D$ is generated by $\lambda_{21}\ldots\lambda_{n1}$,
which is conjugated to $\lambda_{21}$ by the automorphism that
fixes $a_1$ and $a_2$ and sends $a_i$ to $a_2a_i$ for $i>2$. Thus
if $K=mD$, then $\lambda_{21}^m\in K$ and hence all
$\lambda^{m}_{ij}$ are in $K$. Thus
$\phi$ factors through
$Q_{n,m}$.
\QED

If $m=1$, then $Q_{n,1}={\Z}/2$, but in general the groups $Q_{n,m}$
are infinite:

\begin{proposition}  For large $m$, the groups
$Q_{n,m}$ are infinite.
\end{proposition}

\begin{proof} There is a natural homomorphism from  $\autn$ to the
   automorphism group of the free object in the
variety of $n$-generator groups of exponent $m$, which
is written $B_{n,m}$, ``the free Burnside group". The image
of this map contains isomorphic copies of
$B_{n-1,m}$, for example  the subgroup consisting
of the  automorphisms obtained by composing
left Nielsen moves on the first generator:
$a_1\mapsto w a_1,\  a_i\mapsto a_i$ for $
i=2,\dots,n$, where $w$ is a word in the letters
$a_2,\dots,a_n$ and their inverses.  If $m$ is sufficiently
large then this group is infinite \cite{Ad}, \cite{Iv}.

Since the image of each
$\lambda_{ij}^m$ in $\Aut(B_{m,n})$ is trivial, the map
   $\autn\to \Aut(B_{m,n})$ factors
through $Q_{n,m}$.
\QED

\end{document}